\documentclass[10pt,letterpaper]{article}
\usepackage[top=0.85in,left=2.75in,footskip=0.75in]{geometry}

\usepackage{mathtools}
\DeclarePairedDelimiter\abs{\lvert}{\rvert}%

\usepackage{amsmath,amssymb}

\usepackage{changepage}

\usepackage[utf8x]{inputenc}

\usepackage{textcomp,marvosym}

\usepackage{cite}

\usepackage{nameref,hyperref}

\usepackage[right]{lineno}

\usepackage{microtype}
\DisableLigatures[f]{encoding = *, family = * }

\usepackage[table]{xcolor}

\usepackage{array}

\newcolumntype{+}{!{\vrule width 2pt}}

\newlength\savedwidth

\raggedright
\usepackage[aboveskip=1pt,labelfont=bf,labelsep=period,justification=raggedright,singlelinecheck=off]{caption}

\makeatletter
\renewcommand{\@biblabel}[1]{\quad#1.}
\makeatother

\date{}

\usepackage{lastpage,fancyhdr,graphicx}
\usepackage{epstopdf}
\fancyheadoffset[L]{2.25in}
\fancyfootoffset[L]{2.25in}

\begin{document}
\vspace*{0.2in}

\begin{flushleft}
{\LARGE
\textbf\newline{Distance Metrics for Gamma Distributions} 
}
\newline
\\
Colin M. McCrimmon\textsuperscript{1*}
\\
\bigskip
\textbf{1} Medical Scientist Training Program, University of California, Irvine, CA, USA
\\
\bigskip

%
%





* cmccrimm@uci.edu

\end{flushleft}

\section*{Abstract}
Here I present the analytic form of two common distance metrics, the symmetrised Kullback-Leibler Divergence and the Kolmogorov-Smirnov statistic, as well as an extension of the Kolmogorov-Smirnov statistic for comparing theoretical gamma distributions. In doing so, I also present the analytic solution to the intersection of two gamma distributions. Lastly, I provide examples that demonstrate the similarity between these distance metrics and their usefulness in describing the separability of gamma distributions.

\section*{The Gamma Distribution}
The general form of a gamma distribution is

\begin{equation} \label{eq:1:1}
f(x) = \frac {1}{\Gamma (k){\theta}^{k}} {x}^{k-1} {e}^{- \frac {x}{\theta}} \, , \hspace*{0.5cm} \forall x \ge 0
\end{equation}

\noindent{where $k>0$, $\theta>0$ are the shape and scale parameters, respectively. Note that $k$ controls the shape of the distribution, while $\theta$ controls the horizontal stretching/compression. The special case when $k \in \{\, \frac{n}{2} \mid n \in \mathbb{N}, n \ge 1 \,\} $ and $\theta=2$ is known as the chi-squared distribution. The mean and mode (location of the maximum) of any gamma distribution can be easily calculated; the mean is $k\theta$, and the mode is $0$ when $k \le 1$ and $(k-1)\theta$ when $k > 1$.}

\section*{The Kullback-Leibler Divergence}
The Kullback-Leibler (KL) divergence is a commonly-used distance metric that is asymmetric and non-negative. From \cite{Penny2001} and \cite{Bauckhage2014}, it is known that the KL divergence of two gamma distributions, $p(x)$ and $q(x)$, is 

\begin{equation} \label{eq:kl:1}
{D}_{KL}( p(x) || q(x) ) = ({k}_{p} - {k}_{q}) \Psi({k}_{p}) + \log \left( \frac{\Gamma({k}_{q})}{\Gamma({k}_{p})} \right) - {k}_{q} \log \left( \frac{{\theta}_{p}}{{\theta}_{q}} \right) + {k}_{p}\frac{{\theta}_{p} - {\theta}_{q}}{{\theta}_{q}}
\end{equation}

\noindent{where $\Psi$ is the digamma function and $p(x)$ and $q(x)$ have shape and scale parameters ${k}_{p}$, ${ \theta }_{p}$ and ${k}_{q}$, ${ \theta }_{q}$, respectively. Again, the KL divergence in Eq. \ref{eq:kl:1} is asymmetric because ${D}_{KL}( p(x) || q(x) ) \neq {D}_{KL}( q(x) || p(x) )$. However, defining ${D}_{SKL} = {D}_{KL}( p(x) || q(x) ) + {D}_{KL}( q(x) || p(x) )$ produces the following symmetric and non-negative version of the KL divergence for gamma distributions}

\begin{equation} \label{eq:kl:2}
\begin{multlined}
{D}_{SKL} = ({k}_{p} - {k}_{q}) (\Psi({k}_{p}) + \log ({\theta}_{p}) - \Psi({k}_{q}) - \log ({\theta}_{q})) + ({k}_{p}{\theta}_{p} - {k}_{q}{\theta}_{q})\frac{{\theta}_{p} - {\theta}_{q}}{{\theta}_{p}{\theta}_{q}}
\end{multlined}
\end{equation}

\noindent{Another commonly-used divergence metric is the Jensen-Shannon divergence, which is symmetric and bounded by $0$ and $\log(2)$. For two distributions $p(x)$ and $q(x)$, the Jensen-Shannon divergence is defined as
	
\begin{equation} \label{eq:kl:3}
{D}_{JS} = \frac{1}{2} {D}_{KL}( p(x) || m(x) ) + \frac{1}{2} {D}_{KL}( q(x) || m(x) )
\end{equation}

\noindent{where $m(x) = \frac{1}{2} (p(x) + q(x))$. However, if $p(x)$ and $q(x)$ are both gamma distributions, the mixture model $m(x)$ is not easily reducible, and ${D}_{JS}$ has no simple, analytic form.}

\section*{The Kolmogorov-Smirnov Statistic}
Typically, the Kolmogorov-Smirnov (KS) statistic is not explicitly reported, but is instead used in the two-sample KS test to determine whether two empirical cumulative distribution functions (CDFs) were generated from the same underlying distribution. Given empirical CDFs ${F}_{n}(x)$ and ${G}_{m}(x)$, the KS statistic is defined as

\begin{equation} \label{eq:ks:1}
{D}_{KS} = \sup_{x} \left( \abs{{F}_{n}(x) - {G}_{m}(x) } \right)
\end{equation}

\noindent{and is compared in the KS test to $\sqrt{ - \frac {1}{2} \log \left( \frac {\alpha}{2} \right) \left( \frac {n+m}{nm} \right) }$, where $n$ and $m$ are the number of samples in each distribution and $\alpha$ is the desired level of significance. However, the KS statistic itself is an intuitive distance metric to compare the overlap of distributions, and is symmetric and bounded by $0$ and $1$. For two general gamma distributions $p(x)$ and $q(x)$, the KS statistic from Eq. \ref{eq:ks:1} can be defined analytically as
	
\begin{equation} \label{eq:ks:2}
{D}_{KS} = \max_{x \ge 0} \Biggl( \Biggl| \underbrace{ \frac { \gamma ({ k }_{p}, \frac {x}{{ \theta }_{p}})}{ \Gamma ({ k }_{p})}}_{P(x)} - \underbrace{ \frac { \gamma ({ k }_{q}, \frac {x}{{ \theta }_{q}})}{ \Gamma ({ k }_{q})}}_{Q(x)} \Biggr| \Biggr)
\end{equation}

\noindent{where $\gamma$ is the lower incomplete gamma function and ${k}_{p}$, ${\theta}_{p}$ and ${k}_{q}$, ${\theta}_{q}$ are the shape and scale parameters from $p(x)$ and $q(x)$, respectively. Note that $P(x)$ and $Q(x)$ are simply the CDFs of $p(x)$ and $q(x)$, respectively, so ${D}_{KS}$ is clearly bounded by $0$ and $1$ with larger values signifying that the two distributions are highly separable (low overlap). The value of $x$ corresponding to ${D}_{KS}$ is a critical point of Eq. \ref{eq:ks:2}. Thus, if ${x}^{ \ast }$ represents all critical points, then }

\begin{equation} \label{eq:ks:3}
\frac {1}{\Gamma ({ k }_{p})} { \left( \frac {{x}^{ \ast }}{{ \theta }_{p}} \right) }^{{ k }_{p}} {e}^{- \frac {{x}^{ \ast }}{{ \theta }_{p}}} = \frac {1}{\Gamma ({ k }_{q})} { \left( \frac {{x}^{ \ast }}{{ \theta }_{q}} \right) }^{{ k }_{q}} {e}^{- \frac {{x}^{ \ast }}{{ \theta }_{q}}}
\end{equation}

\noindent{This equation can be interpreted as the intersection(s) of the two gamma distributions. Also note that four cases arise when solving for critical points ${x}^{ \ast }$ in Eq. \ref{eq:ks:3}. It is obvious that when ${ k }_{p} = { k }_{q}$ and ${ \theta }_{p} = { \theta }_{q}$, there is no unique solution, and ${D}_{KS} = 0$ for all ${x}^{ \ast }\ge0$. When ${ k }_{p} = { k }_{q}$ and ${ \theta }_{p} \neq { \theta }_{q}$, a single real solution exists at 

\begin{equation} \label{eq:ks:4}
{x}^{ \ast } = \frac {{ k }_{p} { \theta }_{p} { \theta }_{q} \log \left( \frac {{ \theta }_{p}}{{ \theta }_{q}} \right) }{{ \theta }_{p} - { \theta }_{q}}
\end{equation}

\noindent{When ${ k }_{p} \neq { k }_{q}$ and ${ \theta }_{p} = { \theta }_{q}$, a single real solution exists at}

\begin{equation} \label{eq:ks:5}
{x}^{ \ast } = { \theta }_{p} { \left( \frac { \Gamma ({ k }_{p})}{ \Gamma ({ k }_{q})} \right) }^{ \frac {1}{{ k }_{p} - { k }_{q}}}
\end{equation}

\noindent{When ${ k }_{p} \neq { k }_{q}$ and ${ \theta }_{p} \neq { \theta }_{q}$, one or two real solutions to ${x}^{ \ast }$ in Eq. \ref{eq:ks:3} exist, and are given by}

\begin{equation} \label{eq:ks:6}
{x}^{ \ast } = \frac {\textbf {W}(\alpha \beta)}{\alpha}
\end{equation}

\begin{equation}
\alpha = \frac {{ \theta }_{p} - { \theta }_{q}}{({ k }_{p} - { k }_{q}){ \theta }_{p}{ \theta }_{q}} \, , \hspace*{0.5cm} \beta = { \left( \frac {\Gamma ({ k }_{p}) {{ \theta }_{p}}^{{ k }_{p}}}{\Gamma ({ k }_{q}) {{ \theta }_{q}}^{{ k }_{q}}} \right) }^{ \frac {1}{{ k }_{p} - { k }_{q}}}
\nonumber
\end{equation}

\noindent{\textbf{W} in Eq. \ref{eq:ks:6} represents the Lambert-W function. $\textbf{W}(\alpha \beta)$ is single-valued when $\alpha \beta \ge 0$ and double-valued when $\alpha \beta < 0$. Given that ${ k }_{p}$, ${ \theta }_{p}$, ${ k }_{q}$, and ${ \theta }_{q}$ are positive values, simulation suggests that $\alpha \beta > -{e}^{-1}$, and therefore that all values of $\textbf{W}(\alpha \beta)$ are real.

Recall that the critical points, ${x}^{ \ast }$, of Eq. \ref{eq:ks:2} are equivalent to the points of intersection of the two gamma distributions, so ${x}^{\ast} = \{ {x}^{ \ast }_{1} \}$ if $p(x)$ and $q(x)$ intersect once and ${x}^{\ast} = \{ {x}^{ \ast }_{1}, {x}^{ \ast }_{2} \}$ if $p(x)$ and $q(x)$ intersect twice. In the first case, ${x}^{ \ast }_{1}$ is the solution for ${D}_{KS}$, while in the second case, either ${x}^{ \ast }_{1}$ or ${x}^{ \ast }_{2}$ or both are solutions for ${D}_{KS}$. However, ${D}_{KS}$ in the latter underestimates the separability of the two gamma distributions since it only utilizes ${x}^{ \ast }_{1}$ or ${x}^{ \ast }_{2}$ but not both simultaneously. Therefore, I suggest an extension of the standard KS statistic defined as

\begin{equation} \label{eq:ks:7}
{D}_{EKS} = \sum_{i=1}^{n} \abs{ \left( P({x}^{ \ast }_{i}) - Q({x}^{ \ast }_{i}) \right) }
\end{equation}

\noindent{where $n$ is the number of times $p(x)$ and $q(x)$ intersect. Note that ${D}_{EKS} = {D}_{KS}$ for $n = 1$ and ${D}_{EKS} > {D}_{KS}$ for $n = 2$.}

\section*{Examples}
In the previous sections, I defined the symmetrised KL divergence (${D}_{SKL}$) and an extension of the standard KS statistic (${D}_{EKS}$) for gamma distributions. To compare both distance metrics visually, five plots (Fig. \ref{fig:1}) were created using values of ${ k }_{p}$, ${ \theta }_{p}$, ${ k }_{q}$, and ${ \theta }_{q}$ arbitrarily chosen from the interval [1,20]. While both metrics are symmetric and non-negative, ${D}_{EKS}$ may be more practical and intuitive because, unlike ${D}_{SKL}$, it is bounded from above.

Note that in the top plot of Fig. \ref{fig:1} $\text{mode} \left( p(x) \right) < \text{mode} \left( q(x) \right)$ and both distributions have a single point of intersection from Eqs. \ref{eq:ks:4}-\ref{eq:ks:6}. In cases such as this, the complement of ${D}_{EKS}$ and ${D}_{KS}$ is simply the sum of the area under $p(x)$ to the left of ${x}^{\ast}_{1}$ and the area under $q(x)$ to the right of ${x}^{\ast}_{1}$. The second plot of Fig. \ref{fig:1} depicts one example where $p(x)$ and $q(x)$ intersect twice. Here, the value of ${D}_{EKS}$ is $0.486$, whereas the value of the standard KS statistic ${D}_{KS}$ is only $0.368$. Again, because ${D}_{EKS}$ makes use of both ${x}^{ \ast }_{1}$ and ${x}^{ \ast }_{2}$, it may better reflect differences between two distributions.

\begin{figure}[!h]
	\begin{center}
		\includegraphics[width=\linewidth]{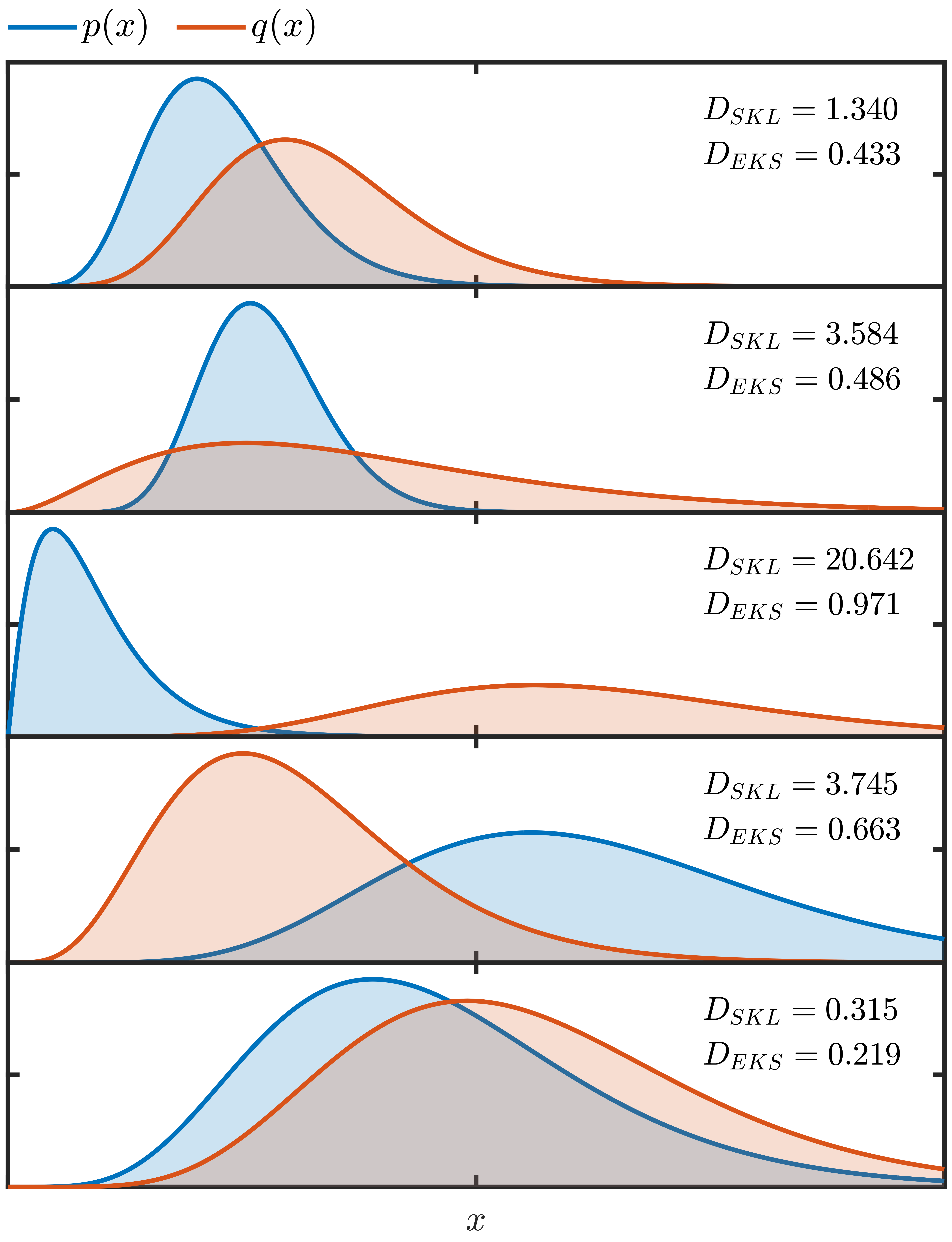}
	\end{center}
	\caption{{\bf Comparison of the symmetric KL divergence (${D}_{SKL}$) and an extension of the KS statistic (${D}_{EKS}$) for gamma distributions.}
		For each plot, arbitrary gamma distributions $p(x)$ (blue) and $q(x)$ (orange) are shown. ${D}_{SKL}$ and ${D}_{EKS}$ were calculated to determine the separability of $p(x)$ and $q(x)$. Note that the y-axis scales are different for each plot, since the total area under each gamma distribution here is unity.}
	\label{fig:1}
\end{figure}


\end{document}